\documentstyle[12pt,amssymb]{amsart}
\newtheorem{lem}{Lemma}[section]
\newtheorem{thm}{Theorem}
\newtheorem{rem}{Remark} 
\newenvironment{proof}
{{\noindent\bf Proof:~}}{\hfill$\Box$}

\def\DD{{\Bbb D}}
\def\RR{{\Bbb R}}
\def\diam{{\mathrm{diam}}}
\def\fr{\frac}
\def\br#1{\left(#1\right)}
\def\brs#1{\left\{#1\right\}}
\def\abs#1{\left|#1\right|}
\def\norm#1{\left\|#1\right\|}
\def\const{{\mathrm{const}}}

\begin{document}

\title{Dimension of quasicircles}
\author{Stanislav Smirnov}
\address{Section de Math\'ematiques, Universit\'e de Gen\`eve\\
2-4 rue du Li\`evre, Case postale 64\\ 
1211 Gen\`eve 4\\ 
Suisse}
\email{Stanislav.Smirnov@@math.unige.ch}
\begin{abstract}
We introduce canonical antisymmetric quasiconformal maps,
which minimize the quasiconformality constant among  maps 
sending the unit circle to a given quasicircle.
As an application we prove Astala's conjecture that
the Hausdorff dimension of a $k$-quasicircle is at most $1+k^2$.
\end{abstract}

\maketitle

A homeomorphism $\phi$ of planar domains is called 
{\em $k$-quasiconformal}, if
it belongs locally to the Sobolev class $W^1_2$
and its {\em Beltrami coefficient}
$$\mu_\phi(z):=\bar\partial\phi(z)/\partial\phi(z),$$
has bounded $L_\infty$ norm: $\norm{\phi}\le k<1$.
Equivalently one can demand that
for almost every point $z$ in the domain of definition
directional derivatives satisfy
$$\max_\alpha\abs{\partial_\alpha f(z)}\le K \min_\alpha\abs{\partial_\alpha f(z)},$$
where the constants of quasiconformality are related by
$$k=\frac{K-1}{K+1}\in[0,1[,~K=\frac{1+k}{1-k}\in[1,\infty[.$$
Quasiconformal maps change eccentricities of infinitesimal ellipses
at most by a factor of $K$,
and it is common to visualize a quasiconformal map $\phi$ 
by considering a measurable field $M(z)$
of infinitesimal ellipses, which is mapped by $\phi$
to the field of infinitesimal circles.
More rigorously, $M(z)$ is an ellipse on the tangent bundle, centered at $0$
and defined up to homothety.
Being preimage of a circle under the differential of $\phi$ at $z$,
in the complex coordinate $v\in T_z{\Bbb C}$
it is given by the equation
$\abs{v+\bar v \mu_\phi(z)}=\const$,
and so its eccentricity is equal to $\abs{1+\abs{\mu_\phi(z)}}/\abs{1-\abs{\mu_\phi(z)}}\le K$.

Quasiconformal maps constitute an important generalization of conformal maps,
which one obtains for $K=1$ (or $k=0$).
However, while conformal maps preserve the Hausdorff dimension,
quasiconformal maps can alter it.
Understanding this phenomenon was a major challenge
until the work \cite{astala-area} of Astala,
where he obtained sharp estimates for the area and dimension distortion 
in terms of the quasiconformality constants.
Particularly,
the image of a set of Hausdorff dimension $1$
under $k$-quasiconformal map
can have dimension at most $1+k$,
which can be attained for certain Cantor-type sets
(like the Garnett-Ivanov square fractal).

Nevertheless this work left open the question of
the dimensional distortion of the subsets of smooth curves.
Astala's \cite{astala-area} implies
that a $k$-quasicircle (that is an image of a circle
under a $k$-quasiconformal map) has the Hausdorff dimension
at most $1+k$, being an image of a $1$-dimensional set.
However, some years earlier,
Becker and Pommerenke have shown in \cite{becker-pommerenke}
that the Hausdorff dimension of a $k$-quasicircle
cannot exceed $1+37k^2$. 
The Becker-Pommerenke estimate is better than Astala's
for very small values of $k$ because of its quadratic behavior,
which they have shown to be sharp.
Apparently, some properties of a circle were making
it harder to increase its dimension by a quasiconformal map.
Motivated by the two available estimates,
Astala conjectured in \cite{astala-area}
that the Hausdorff dimension of $k$-quasicircles cannot
exceed $1+k^2$.

In this paper we prove Astala's conjecture:
\begin{thm}\label{thm:astala}
The Hausdorff dimension of a $k$-quasicircle
is at most $1+k^2$.
\end{thm}
\begin{rem}
Our proof has some technical simplifications compared to Astala's since we are estimating only the dimension,
and do not address the question of the Hausdorff measure.
\end{rem}
\begin{rem}
It easily follows that a $k$-quasiconformal image of any smooth curve
satisfies the same estimate.
\end{rem}

The question  about the sharpness remains open.
In contrast to the general case, where sharpness of
Astala's estimates follows relatively easily from the proof,
here the question seems to be more difficult
and have deep connections.
Our methods rely on representing quasicircles by
canonical antisymmetric Beltrami coefficients,
and are related to difficult questions
of quantifying fine structure of harmonic measures
and boundary dimensional distortion of conformal maps.

The strategy of our proof is similar to Astala's \cite{astala-area}:
we employ holomorphic motions and
thermodynamic formalism of \cite{ruelle-book,ruelle-repellers}.
The difference is that we establish a ``canonical Beltrami representation''
for quasicircles, which exhibits an interesting kind of symmetry.
This allows us to embed the quasicircle into a symmetric holomorphic
motion and use a better version of Harnack's inequality.

We formulate the Theorem for a $k$-quasiline,
that is an image of the real line ${\Bbb R}$ under a $k$-quasiconformal map.
One obtains the quasicircle formulation by conjugating by a M\"obius map
which sends the real line to a circle and hence replaces the complex conjugation $\bar z$
by the inversion with respect to the circle in question.

\begin{thm}[Canonical Beltrami representation]\label{thm:repr}
For a curve $\Gamma$ the following are equivalent:
\renewcommand{\theenumi}{{\roman{enumi}}}
\begin{enumerate}
\item \label{i} $\Gamma$ is a $k$-quasiline.
\item \label{ii} $\Gamma=\psi({\Bbb R})$ with $\|\mu_\psi\|\le 2k/(1+k^2)$
and $\mu_\psi(z)=0,~z\in{\Bbb C}_{+}$.
\item \label{iii} $\Gamma=\phi({\Bbb R})$ with $\|\mu_\phi\|\le k$
and 
\begin{equation}\label{eq:asymm}
\mu_\phi(\bar z)=-\overline{\mu_\phi(z)}~.
\end{equation}
\end{enumerate}
\end{thm}

\begin{rem}
In the ellipses language condition (\ref{eq:asymm}) means that 
the ellipse $M(\bar z)$ at $z$ is a rotated by $\pi/2$ 
complex conjugate of the ellipse $M(z)$. 
\end{rem}

\begin{rem}\label{rem:symm}
If 
\begin{equation}\label{eq:symm}
\mu_\phi(\bar z)=\overline{\mu_\phi(z)}~,
\end{equation}
or equivalently the infinitesimal ellipse $M(z)$ at $z$
is a complex conjugate $\overline{M(\bar z)}$ of the ellipse at $\bar z$,
then by symmetry the real line is mapped to the real line
(though quasisymmetrically perturbed),
so Beltrami coefficient is ``wasted.''
On the other hand, the Theorem above shows
that in case of (\ref{eq:asymm})
Beltrami coefficient is ``optimal''
among ones representing the same quasicircle.
Note that we represent every quasiconformal map
as a superposition of two maps with
better quasiconformality constants
and satisfying (\ref{eq:asymm}) and (\ref{eq:symm}) respectively.
\end{rem}

\paragraph{\bf Acknowledgments}
This research was partially supported by the NSF grant DMS-9706875
as well as by the EU RTN CODY and Swiss NSF.

\section{Proofs}

\subsection{Proof of the Theorem~\ref{thm:repr}}\label{sec:repr}~  

It is more illustrative to work with ellipse fields rather than 
Beltrami coefficients.
For a (measurable) ellipse field $M(z)$ we will denote by
$\|M\|$ the essential supremum of the eccentricities,
which is related to the norm of the corresponding Beltrami coefficient
$\mu$ by $\|M\|=(1+\|\mu\|)(1-\|\mu\|)$.

\noindent{$\bf(\ref{i})\Rightarrow(\ref{ii})$.}
Let $N(z)$, $\|N\|\le K:=(1+k)/(1-k)$, be the ellipse field
representing a $k$-quasiconformal map $\eta$,
which maps ${\Bbb R}$ onto $\Gamma$. 
We will ``correct'' the map $\theta$ by precomposing
it with a carefully selected map which preserves the real line.
We construct this new quasiconformal map $\alpha$
by the means of an ellipse field $A$:
$$
A(z)=\left\{
\begin{matrix}
\overline{N(\bar z)},~z\in{\Bbb C}_{-} \\
N(z),~z\in{\Bbb C}_{+}\\
\end{matrix}\right.~.
$$
Then the map $\psi:=\eta\circ\alpha^{-1}$ is the desired one.
Indeed, by the Remark~\ref{rem:symm} the map $\alpha$ preserves the real line,
so $\psi({\Bbb R})=\eta(\alpha^{-1}({\Bbb R}))=\eta({\Bbb R})=\Gamma$.
For $z$ in the upper half plane both $\eta$ and $\alpha$ send
the ellipse field $N(z)$ to the field of circles, hence
the map $\psi=\eta\circ\alpha^{-1}$ preserves the field of circles
and is conformal in the upper half-plane.
In the lower half plane both $\eta$ and $\alpha$
change eccentrivities by at most $K$, so
$\psi$ changes eccentricities by at most $K^2$.
Thus $\|\mu_\psi\|\le (K^2-1)/(K^2+1)=2k/(1+k^2)$.

\noindent{$\bf(\ref{ii})\Rightarrow(\ref{iii})$.}
Let $M(z)$, $\|M\|\le K'$, be the ellipse field
representing the $2k/(1+k^2)$-quasiconformal map $\psi$,
with
$$K'~=~\frac{1+2k/(1+k^2)}{1-2k/(1+k^2)}~=~
\br{\frac{1+k}{1-k}}^2~=~K^2~.$$
Define a new quasiconformal map $\beta$
by the means of a new ellipse field $B$:
$$
B(z)=\left\{
\begin{matrix}
\sqrt{M(z)},~z\in{\Bbb C}_{-} \\
\sqrt{\overline{M(\bar z)}},~z\in{\Bbb C}_{+}\\
\end{matrix}
\right.~.
$$
Here $\sqrt M$ denotes the ellipse 
with the same alignment whose eccentricity
is the square root of $M$'s eccentricity
(note that we do not rotate the ellipses, just change their eccentricities).
Then the map $\phi:=\psi\circ\beta^{-1}$ is the desired one.
As before, by the Remark~\ref{rem:symm} the map $\beta$ preserves the real line,
so $\phi({\Bbb R})=\psi({\Bbb R})=\Gamma$.

Let $L(z)$ be the image of the field $M(z)$ under the map $\beta$,
or equivalently the preimage of the circle field
under the map $\phi$.
For $z$ in the lower half-plane  $\beta$ sends
the ellipse field $\sqrt{M(z)}$ to the field of circles, hence
we conclude that there $\|L\|\le K=\sqrt{K'}$.
Using ``symmetric'' construction of $\beta$ we conclude that
for $z$ in the upper half-plane, it sends
$M(z)$ which is just a circle field there
to $i L(z)$, i.e. the rotation of $L(z)$ by $\pi/2$. 
Rephrasing in terms of Beltrami coefficients, we
arrive at (\ref{eq:asymm}).

\noindent{$\bf(\ref{iii})\Rightarrow(\ref{i})$.}
Trivial since $\phi$ is already $k$-quasiconformal.
\hfill$\Box$

\subsection{Symmetric Harnack}\label{sec:harn}

\begin{lem}[Symmetric Harnack]\label{lem:harnack}
Let $h$ be a positive harmonic function in the unit disc ${\Bbb D}$,
whose partial derivative at the origin vanishes in the direction
of some $\lambda\in{\Bbb D}$: $\partial_\lambda h(0)=0$.
Then $h$ satisfies the following stronger
version of Harnack's inequality:
$$(1-|\lambda|^2)/(1+|\lambda|^2)h(\lambda)\le
h(0)\le(1+|\lambda|^2)/(1-|\lambda|^2)\,h(\lambda)~.$$
\end{lem}

\begin{rem}
Unlike the usual Harnack, the one-sided estimate does not
remain valid for subharmonic functions. 
Therefore the combination of Becker-Pommerenke and Astala's results
does not automatically give the $1+k^2$ estimate.
\end{rem}

\begin{proof} 
Without loss of generality we can assume that $h(0)=1$
and gradient of $h$ vanishes at the origin.
Otherwise we consider the function $(h(z)+h(z^*))/(2 h(0))$,
where ${}^*$ denotes the symmetry with respect to the line
through the origin and $\lambda$.

Then the function $f:=\frac{z-1}{z+1}\circ(h+i\tilde h)$
is analytic from the unit disc ${\Bbb D}$ to itself
(here $\tilde{h}$ denotes the harmonic conjugate).
Moreover, $f$ preserves the origin and $f'(0)=0$,
since$\nabla h(0)=0$.
Thus we can apply the Schwarz lemma to the function $f(z)/z$,
concluding that $\abs{f(\lambda)}\le\abs{\lambda}^2$.
This implies that
$h(\lambda)\le(1+\abs{\lambda}^2)/(1-\abs{\lambda}^2)$
and
$h(\lambda)\ge(1-\abs{\lambda}^2)/(1+\abs{\lambda}^2)$.
\end{proof}

\subsection{Proof of the Theorem~\ref{thm:astala}}\label{sec:dim}~  

Take a $k$-quasiline $\Gamma$ and
its representation by a $k$-quasiconformal map $\phi$
with Beltrami coefficient $\mu$ satisfying (\ref{eq:asymm}).
Include $\phi$ into the standard holomorphic motion
of maps $\phi_\lambda,\,\lambda\in{\Bbb D}$,
whose Beltrami coefficients are given by
$\mu_\lambda := \mu \cdot \lambda/k$
and which preserve points $0$, $1$, $\infty$.
As usual, $\phi_0={\mathrm{id}}$ and $\phi_k=\phi$.
Because of the particular symmetric behavior (\ref{eq:asymm}) of
$\mu_k$ we conclude that
$\mu_\lambda$ satisfies (\ref{eq:asymm}) for real $\lambda$
and (\ref{eq:symm}) for imaginary $\lambda$.
There is additional symmetry:
for real values of $\lambda$ one has $\phi_\lambda(z)~=~\overline{\phi_{-\lambda}(\bar z)}$,
which follows from the same property for $\mu_{\lambda}$. 

We proceed as in \cite{astala-area} with some simplifications in estimating the distortion.

To this effect fix a radius $\rho\in(1/2,1)$, and consider for the time being
only $\lambda$ inside the slightly smaller disk $\rho\DD$,
so the maps $\phi_\lambda$ have
uniform quasiconformality constants.
Note that by the results of Gehring (see e.g. Theorem 18.1 in \cite{vaisala})
the maps $\phi_\lambda$ are uniformly quasisymmetric, so there is a constant $C=C(\rho)$
such that 
\begin{eqnarray}\label{eq:qs}
\abs{z-x}\le\abs{y-x}~\Rightarrow~\abs{\phi_\lambda(z)-\phi_\lambda(x)}
\le C\cdot\abs{\phi_\lambda(y)-\phi_\lambda(x)},\\
\label{eq:qss}
C\cdot \abs{z-x}\le\abs{y-x}~\Rightarrow~2 \abs{\phi_\lambda(z)-\phi_\lambda(x)}
\le \abs{\phi_\lambda(y)-\phi_\lambda(x)}.
\end{eqnarray}

It is enough to estimate the dimension of the image of the interval $[0,1]$.
Cover it by $n$ length $1/n$ intervals $[a_j,b_j]$ and denote
by $B_j(\lambda)$ the ball centered at $\phi_\lambda(a_j)$
whose boundary circle passes through $\phi_\lambda(b_j)$.
Note that its ``complex radius''
$$r_j(\lambda):=\phi_\lambda(b_j)-\phi_\lambda(a_j),$$
is a holomorphic function of $\lambda$.

The image of the interval $[0,1]$ is covered 
by the images of small intervals, 
each having diameter at most $C \abs{r_j(\lambda)}$ by (\ref{eq:qs}).
So to estimate the Hausdorff dimension we have to estimate the sums
\begin{equation}\label{eq:dist}
\sum_j \diam\br{\phi_\lambda[a_j,b_j]}^p
\le C^p\,\sum_j \abs{r_j(\lambda)}^p.
\end{equation}
We will estimate the logarithm of the right-hand side.
The logarithm is a concave function, so applying the Jensen inequality 
to the collection of points $\brs{\abs{r_j(\lambda)}^p/\nu_j}$ with weights $\brs{\nu_j}$ summing to one, 
we obtain
\begin{equation}\label{eq:jensen}
\log \sum_j \abs{r_j(\lambda)}^p
= \log \sum_j  \nu_j \fr{\abs{r_j(\lambda)}^p}{\nu_j} 
\ge \sum_j \nu_j \log{\fr{\abs{r_j(\lambda)}^p}{\nu_j}}
={I_\nu-p\Lambda_\nu(\lambda)}~,
\end{equation}
where $I_\nu:=-\sum_j \nu_j \log \nu_j$
is the ``entropy'' and
$\Lambda_\nu(\lambda):=-\sum_j  \nu_j \log \abs{r_j(\lambda)}$
is the ``Lyapunov exponent'' of the probability
distribution $\brs{\nu_j}$.
Note that the Lyapunov exponent is a harmonic function of $\lambda$,
since $r_j(\lambda)$ are holomorphic.

The equality in (\ref{eq:jensen}) is attained if $\nu_j$ are chosen
proportionalto $\abs{r_j(\lambda)}^p/\nu_j$ coincide,
and we deduce
the {\em variational principle} from Astala's \cite{astala-area}:
\begin{equation}\label{eq:varpr}
\log \sum_j \abs{r_j(\lambda)}^p
=\sup_{\nu}\brs{I_\nu-p\Lambda_\nu(\lambda)}~,
\end{equation}
where the supremum is taken over all probability distributions $\nu$.

Fix some distribution $\nu$ and consider the function
$$H(\lambda):=2\Lambda_\nu(\lambda)-I_\nu+3\log C.$$
It is harmonic in $\lambda$ (since $\Lambda_\nu$ is)
and is an even function on the real line (because of the symmetry of our motion
$r_j(\lambda)=\overline{r_j(-\lambda)}$ for $\lambda\in\RR$).
By (\ref{eq:qss}) the balls $B_j(\lambda)$ cover every point at most $C$ times,
while by (\ref{eq:qs}) their  union is contained in a ball of radius $C$.
Hence we can write
$$
\sum_j\abs{r_j(\lambda)}^2\le C^3,
$$
deducing by the variational principle (\ref{eq:varpr}) that
$$
I_\nu-2\Lambda_\nu(\lambda)\le\log C^3,
$$
and concluding that $H$ is non-negative in the disk ${\rho}\DD$.
Finally,
$$I_\nu-\Lambda_\nu(0)\le\log\sum{\abs{r_j(0)}} = \log1=0,$$
therefore
$$H(0)=2\Lambda_\nu(0)-I_\nu+3\log C\ge I_\nu+3\log C.$$

We can now apply the Lemma~\ref{lem:harnack} (in the smaller disk $\rho\DD$)
to obtain
$$2\Lambda_\nu(k)-I_\nu+3\log C=H(k)\ge
\frac{1-k^2\rho^{-2}}{1+k^2\rho^{-2}}H(0)
\ge \frac{1-k^2\rho^{-2}}{1+k^2\rho^{-2}}\br{I_\nu+3\log C},$$
and deduce
$$\frac{2}{1+k^2\rho^{-2}}I_\nu-2\Lambda_\nu(k)\le
\frac{2 k^2\rho^{-2}}{1+k^2\rho^{-2}}{3\log C},$$
which can be rewritten as
$${I_\nu-\br{1+k^2\rho^{-2}}\Lambda_\nu(k)}
\le{k^2\rho^{-2}}{3\log C}.$$
Since the latter holds for all distributions $\nu$,
we conclude by the variational principle (\ref{eq:varpr}) that
\begin{equation}
\log\sum_j\abs{r_j(k)}^p
=\sup_\nu\brs{I_\nu-p\Lambda_\nu(k)}
\le{k^2\rho^{-2}}{3\log C}\le 12 \log C,
\label{eq:est}
\end{equation}
where we set $p:=1+k^2\rho^{-2}$.

Sending $n$ to infinity, we conclude by (\ref{eq:dist}) and (\ref{eq:est}) that the $p$-dimensional
Hausdorff measure of the quasiline $\phi[0,1]$ is bounded by $C^{14}$,
and hence its dimension is at most $p=1+k^2\rho^{-2}$. 
Sending now $\rho$ to $1$ we obtain the desired estimate.
\hfill$\Box$

\section{Quasisymmetric distortion}\label{sec:qsym}  
For images of subsets of the circle the methods above cannot be applied to the full
extent (since symmetrizing Beltrami preserves image of the circle but not of 
particular subsets). 

But if the circle is preserved (e.g. we are dealing with
maps which restrict to quasisymmetric maps of the unit circle),
we can work with the same symmetrized holomorphic motions,
only we are exploring a direction orthogonal to the previous one.
We say that a mapping of the unit circle is {\it $k$-quasisymmetric},
if it can be extended to a $k$-quasiconformal map of the complex plane.
Note that this quasisymmetry constant is different from the usual $\rho$.

In this case, if we start with a subset of the circle of dimension $1$, its dimension can never be bigger 
than one when $\lambda$ takes real values, so pressure-type quantities have zero gradient,
and the same symmetric Harnack gives that $k$-quasisymmetric image has dimension at least $1-k^2$
(which is the ``dual'' statement to the one that in the orthogonal direction we get at most $1+k^2$).
This and related results are discussed by Prause in \cite{prause}.

If we start with a subset of the circle 
of dimension $a<1$, we can use that its dimension will be 
at most $1$ for real $\lambda$, 
and improve Harnack's inequality using this global information.
The corresponding estimates are rather cumbersome, and will be the subject of a separate paper.

\end{document}